\documentclass[11pt, dvipsnames]{article}
\usepackage{amsmath,amsthm,amssymb,mathtools,amsfonts}
\usepackage{mathrsfs}
\usepackage{hyperref}
\usepackage[capitalise]{cleveref}
\usepackage{color}
\usepackage[table,xcdraw]{xcolor}
\usepackage{geometry}
\usepackage{placeins}
\usepackage{booktabs}
\usepackage{caption}
\usepackage{subcaption}

\usepackage{algorithm}

\usepackage{multicol}
\usepackage{algpseudocode}

\newcommand{\act}{\mathrm{act}}

\newcommand{\T}{\top}

\usepackage{subcaption}

\usepackage{accents}

\newcommand{\roundup}[1]{\overline{#1}}
\newcommand{\rounddown}[1]{\underline{#1}}

\newcommand{\miplib}{\mbox{MIPLIB}}
\Crefname{theorem}{Theorem}{Theorems}
\Crefname{proposition}{Proposition}{Propositions}
\Crefname{lemma}{Lemma}{Lemmas}
\Crefname{section}{Section}{Sections}
\Crefname{appendix}{Section}{Sections}
\Crefname{definition}{Defition}{Definitions}

\usepackage{xspace}

\newcommand{\solver}[1]{\textsc{#1}\xspace}

\newcommand{\scipversion}{8.0.3}

\newcommand{\scip}{\solver{SCIP}}
\newcommand{\scipv}{\solver{SCIP}~\scipversion\xspace}

\newcommand{\soplexversion}{6.0.3}
\newcommand{\soplexv}{\solver{SoPlex}~\soplexversion\xspace}

\newcommand{\papiloversion}{2.0.1}
\newcommand{\papilov}{\solver{PaPILO}~\papiloversion\xspace}

\usepackage[utf8]{inputenc}
\usepackage[english]{babel}
\usepackage{float}
\usepackage[affil-it]{authblk}

\title{Certified Constraint Propagation and Dual Proof Analysis in a Numerically Exact MIP Solver}
\author[1]{Sander Borst}
\author[2]{Leon Eifler}
\author[2,3]{Ambros Gleixner}

\affil[1]{Centrum Wiskunde \textsl{\symbol{`\&}} Informatica, The Netherlands, \url{sander.borst@cwi.nl}}
\affil[2]{Zuse Institute Berlin, Germany, \url{eifler@zib.de}}
\affil[3]{HTW Berlin, Germany, \url{gleixner@htw-berlin.de}}

\begin{document}
\maketitle

\begin{abstract}
	This paper presents the integration of constraint propagation and dual proof analysis in an exact, roundoff-error-free MIP solver. The authors employ safe rounding methods to ensure that all results remain provably correct, while sacrificing as little computational performance as possible in comparison to a pure floating-point implementation. The study also addresses the adaptation of certification techniques for correctness verification. Computational studies demonstrate the effectiveness of these techniques, showcasing a 23\% performance improvement on the \miplib~2017 benchmark test set.
\end{abstract}
\section{Introduction}

Mixed integer programming (MIP) is a powerful modeling technique that can be used to solve a wide variety of complex problems.
Algorithms that solve MIP problems, short MIPs, use an LP-based branch and bound algorithm at their core, but employ many complementary techniques to improve their performance~\cite{NemhauserWolsey1988,Achterberg2007,ConfortiCornuejolsZambelli2014,AchterbergBixby2019,HojnyPfetsch2019}.
They typically use double precision floating-point arithmetic together with the careful handling of numerical tolerances to quickly compute accurate solutions.
While this approach makes the solving process highly efficient, and is accurate enough in practice for most applications, there exist problems where exact solutions are necessary.
This is the case for applications in computational mathematics~\cite{Bofi19,Burt12,EiflerGleixnerPulaj2022,KenterEtAl2018,LanciaEtAl2020,Pulaj20} as well as for several industry applications where exact correctness is critical~\cite{Achterberg2007,WilkenEtAl2000,SahraouiBendottiAmbrosio2019}.
Furthermore, there exist pathological instances that exhibit such numerical difficulties that floating-point solvers produce large errors.
For these reasons, there is a need for roundoff-error-free MIP solvers.

In recent years, a hybrid approach has been proposed that aims to take advantage of both floating-point and exact arithmetic in order to solve MIPs exactly \cite{cook_hybrid_2013}.
This approach has been revised and further improved by \cite{eifler_computational_2022}, and more recently in~\cite{EiflerGleixner2023} with the addition of cutting planes.
What can be clearly seen in these works is that for numerical techniques to be beneficial, it is crucial to employ them as much as possible using floating-point arithmetic in a numerically safe way, as symbolic computation can be prohibitively slow.

An important feature of MIP solvers is \emph{constraint propagation}, which is a technique that uses the problem constraints to tighten variable bounds, which can be used to detect infeasibility of a subproblem before the LP is solved. This saves time, as performing constraint propagation is generally much faster than solving an LP.

Another key feature of MIP solvers that has not yet been adapted to the exact setting is \emph{conflict analysis}, which can be categorized in two distinct types.
In this paper we will study \emph{dual proof analysis}, which was introduced by \cite{witzig_experiments_2016,witzig_computational_2021}.
This type of conflict analysis derives constraints from Farkas certificates of infeasible subproblems and dual solutions of bound-exceeding LPs. These are redundant constraints that do not strengthen the LP relaxation. However, by performing propagation on these constraints, infeasibility of subsequent subproblems can be detected more efficiently.
Another form of conflict analysis is \emph{graph-based} conflict analysis \cite{achterberg_conflict_2007}, which was inspired by a similar procedure in SAT solving \cite{jr_using_1997, marques-silva_grasp_1999, moskewicz_chaff_2001}.
This form of conflict analysis is numerically simpler, as it only generates disjunctive constraints with $\pm 1$ coefficients.
In this paper we study dual proof analysis, since our focus is on adapting numerical methods and studying the challenges that arise in the exact setting.

An important aspect in the context of roundoff-error-free algorithms is the \emph{certification of their correctness}.
A certifying algorithm creates a certificate alongside its solution that can be \emph{independently verified} to be correct~\cite{McConnellMehlhornEtAl2011,AlkassarEtAl2011}.
Such techniques have become standard in SAT solving~\cite{HeuleHuntWetzler2014,Cruz-FilipeHeuleEtAl2017,Cruz-FilipeMarquesEtAl2017,Goldberg2003}, and have also been adapted to SMT solving~\cite{de2008proofs,BarbosaEtal2022}, pseudo boolean optimization~\cite{Gocht_2019_veripb,ElffersGochtMcCreeshNordstrom_2020}, as well as to MIP solving~\cite{cheung_verifying_2017,eifler_computational_2022,EiflerGleixner2023}.

To summarize, our contribution is to develop numerically safe versions constraint propagation and dual proof analysis, and to show that they can be used to improve the performance of an exact MIP solver in practice. Furthermore, we show how to certify the correctness of the derived bounds in the VIPR~\cite{cheung_verifying_2017} certificate format.
Our paper is structured as follows. First, we introduce the existing methods of dual proof analysis and propagation in \cref{sec:review}.
Then, we describe how to adapt these techniques to the exact setting as well as how to certify their correctness in \cref{sec:exact}.
We conduct a thorough computational study in \cref{sec:comp}, showing that these techniques can be used to improve the running time by 23\% on the \miplib~2017 benchmark test set~\cite{MIPLIB} and solve more instances within a time limit of two hours.
We conclude with an outlook on future research in \cref{sec:conclusion}.

\section{Constraint Propagation and Conflict Analysis}
\label{sec:review}
In this paper we consider MIPs of the form
\begin{align*}
	\text{minimize}\hspace{5.2em} c^\T x            &                         \\
	\text{subject to }\hspace{3em} Ax\geq b&,                         \\
	\ell_i\leq x_i\leq u_i & \quad \text{for all } i \in [n], \\
	x_i\in \mathbb{Z}      & \quad \text{for all } i\in S,
\end{align*}
for $A\in \mathbb{R}^{m\times n}, 	c\in \mathbb{R}^{n},b\in \mathbb{R}^{m}$. For the sake of simplicity of presentation, we assume that all variables have finite bounds. The techniques are still applicable to MIPs with infinite bounds, but propagation of some constraints involving unbounded variables may become ineffective and is skipped.
In the following, we first review the techniques of constraint propagation and dual proof analysis in the floating-point setting.

\subsection{Constraint Propagation}
\label{sec:propagationfp}

Constraint propagation is a technique that uses the problem constraints to tighten variable bounds, see, e.g. \cite[p.~426]{Achterberg2007}.
Consider a constraint $a^\T x \leq b$.
The maximum activity is defined as $\act^+:=\sum_{i=1}^n\max(a_i\ell_i, a_iu_i)$. Similarly, the minimum activity is defined to be $\mathrm{act}^-:=\sum_{i=1}^n\min(a_i\ell_i, a_iu_i)$.
We define the maximal activity relative to variable $x_k$ as
\begin{equation*}
	\act_k^+:=\sum_{i\in [n]\setminus \{k\}}\max(a_i\ell_i, a_iu_i),
\end{equation*}
and the minimal activity $\act_k^-$ analogously.

Any feasible vector $x$ will satisfy $a^\T x \leq \act^+$.
Hence, for any variable $x_i$ with $a_k>0$, we have
\begin{align}
	x_k \geq \frac{b-\act_k^+ }{a_k},
	\label{eq:act_bound1}
\end{align}
and if $a_k<0$ we have
\begin{align}
	x_k \leq \frac{b-\act_k^+ }{a_k}.
	\label{eq:act_bound2}
\end{align}

Constraint propagation is the technique of computing the appropriate bounds from \cref{eq:act_bound1,eq:act_bound2} and using them to tighten one of the variable bounds $\ell_k, u_k$ if possible. It can be performed in an iterated fashion, as long as some of the bounds have been tightened. This is also the case at every node in a branch-and-bound tree, where the local bounds of the subproblem are taken into account after they have been tightened by branching decisions.

In addition, by computing the minimal activity of a constraint, it can sometimes be determined that a subproblem is infeasible. This happens when $\act^-< b$. In this case, the problem is infeasible because the constraint is violated by any solution that respects the variable bounds.
Conversely, whenever $\act^+ \geq b$, then constraint is redundant and can be safely removed for the corresponding subproblem.

\subsection{Dual Proof Analysis}
\label{sec:dualprooffp}

Dual proof analysis is a technique that allows to learn from previously solved subproblems that are infeasible or exceed the objective bound.
In this type of conflict analysis, for any subproblem with an infeasible or bound-exceeding LP relaxation a constraint is added to the global problem from which the infeasibility of the subproblem can be derived by applying constraint propagation to this single constraint instead of solving the LP relaxation over all constraints.
Note that the constraint will be redundant, and therefore does not need to be added to the LP.
However, it can be used for constraint propagation in other nodes of the tree.

If branching is only performed on variables, then such a subproblem is given as a sequence of bounds $\ell'_i, u'_i$ that have to be satisfied in addition to the problem constraints. So, the LP relaxation in any node in the tree will be of the form
\begin{align*}
	\min c^\T x \\Ax\geq b \\\ell_i'\leq x_i\leq u_i' &\quad \text{for all } i\in [n].
\end{align*}
The dual program to this LP is given by
\begin{align*}
	\max\ b^\T y+ \ell'^\T r^+ -u'^\T r^-  \\
	r^+ -r^- - A^\T y = -c                \\
	y\geq 0, r^+\geq 0, r^-\geq 0.
\end{align*}
If  $\ell'\leq u'$, we can assume that for all $i\in [n]$ either $r_i^+$ or $r_i^-$ is zero, so we will write $r=r^+-r^-$.
If the subproblem is infeasible, then the dual program will be unbounded. This implies that there must exist a ray $y,r^+, r^-$ with
\begin{align*}
	b^\T y + \ell'^\T r^+-u'^\T r^-  & > 0, \\
	A^\T y + r                       & = 0.
\end{align*}
Such a ray is called a Farkas proof of infeasibility \cite{farkas_theorie_1902}. Rewriting this gives
\begin{align}
	0 < b^\T y +\ell'^\T r^+-u'^\T r^-  = b^\T y -  \ell'^\T (A^\T y)^- + u'^\T (A^\T y)^+.
	\label{eq:conflict}
\end{align}
In dual proof analysis, the MIP solver adds the constraint $y^\T A x\leq y^\T b$ to the problem.
Note that this constraint is globally valid, since it is a linear combination of problem constraints.
The minimum activity of this constraint in the current subproblem is $\act^-= \ell'^\T (A^\T y)^+-u'^\T (A^\T y)^- $. Substituting this into \cref{eq:conflict} shows that $\act^- > b^\T y$. Hence, for this subproblem, infeasibility can immediately be derived using propagation of the newly added constraint.

Now consider the case of a bound-exceeding subproblem. This occurs when the objective value to the LP relaxation of the current problem is higher than the objective value $c^\T x^\bullet$ of the best IP-solution $x^\bullet$ found so far. We can cut off bound-exceeding subproblems by imposing $c^\T x \leq c^\T x^\bullet$.
Let $y$ be the optimal dual solution. Dual proof analysis adds the constraint $(y^\T A- c)^\T x \geq y^\T b -c^\T x^\bullet$.
This constraint is valid, as it is a linear combination of the constraints of $A$ and $c^\T x \leq c^\T x^\bullet$.
Note that for any feasible $x$ for the subproblem, we have
\begin{align}
	b^\T y+ \ell'^\T (A^\T y -c)^+ -u'^\T (A^\T y -c)^- =  b^\T y +\ell'^\T r^+ -u'^\T r^- \geq c^\T x^\bullet,
\end{align}
that is
\begin{align}
  \ell'^\T (A^\T y -c)^+-u'^\T (A^\T y -c)^- \geq c^\T x^\bullet - b^\T y.
  \end{align}
Note that the left hand is the minimum activity of the newly added conflict constraint. So just like in the case of an infeasible subproblem, we see that the bound-exceeding subproblem can now be cut of using constraint propagation on the newly added constraint.

The constraints will not be added to the LP, but are only used for propagation.
Since dual proof analysis might be applied many times it can get too expensive to store all the derived constraints.
Therefore, we employ a method called aging, to dynamically remove constraints that have not been useful for a long time~\cite{witzig_experiments_2016}.
In the exact setting, we make use of the same method.

\section{Application in Numerically Exact MIP Solving}
\label{sec:exact}
To efficiently implement constraint propagation and dual proof analysis in an exact MIP solver, we will perform most computations in floating-point arithmetic. To guarantee correctness, we carefully make use of directed rounding.
We start with some definitions.

Let $\mathbb{F} \subseteq \mathbb{Q}$ denote the set of floating-point numbers. In practice, these will be standard IEEE double-precision \cite{IEEE754-2008} numbers with 11 bits for the exponent and 52 bits for the mantissa. For all $x\in \mathbb{Q}$, we define:
\begin{align*}
	\roundup{x} = \min\{y\in \mathbb{F}: y \ge x\} && \rounddown{x} = \max\{y\in \mathbb{F}: y \le x\}
\end{align*}
Whenever we write $\roundup{a+b+c+\cdots}$ we mean $\roundup{\roundup{\roundup{a}+b}+\cdots}$.
An expression $\rounddown{a+b+\cdots}$ is defined analogously, as well as multiplication, division, and any combination thereof. We note that this is consistent with how these expressions are computed in practice.

\subsection{Numerically Safe Constraint Propagation}
Constraint propagation is applied at every node in the branch-and-bound tree. Recomputing the minimum and maximum activity each time would be costly. Therefore, a solver will keep track of these activities and update them whenever they change.
This happens when the variable bounds change, i.e.,
when the bounds $(\ell, u)$ of variable $x_i$ become $(\ell', u')$, we update the maximum activity as
\begin{align*}
	\act^+ \gets \begin{cases}
		(u'-u)a_i       & \text{if } a_i\geq 0, \\
		(\ell'-\ell)a_i & \text{otherwise}.
	\end{cases}
\end{align*}

In the exact setting, the activities would ideally be computed in exact arithmetic, to give the tightest possible bounds.
However, updating the activities using symbolic computations can be prohibitively expensive.
So instead our implementation uses floating-point arithmetic.
To ensure that no incorrect bounds are derived, we enforce the maximum activity that we maintain to be at least as large as the actual maximum activity.
So when bounds $(\ell,u)$ are updated to $(\ell', u')$, the maximum activity is updated to
\begin{align*}
	\act^+ \gets \begin{cases}
		\roundup{\act^++u'a_i-ua_i}       & \text{if } a_i\geq 0, \\
		\rounddown{\act^++\ell'a_i-\ell a_i} & \text{otherwise},
	\end{cases}
\end{align*}
using directed rounding.
Similarly, for the minimum activity we maintain only a lower bound on the exact value.

In the floating-point setting the maintained activities need to be recomputed regularly to keep the aggregated inaccuracy from accumulating. In the exact setting this is not necessary, since the activities are guaranteed to remain valid. For this reason we do not recompute them.

\subsection{Numerically Safe Dual Proof Analysis}
As explained in \cref{sec:dualprooffp}, the dual ray $y$ is computed and the constraint $y^\T A x\leq y^\T b$ is added to the problem. For the sake of efficiency, we let the MIP solver compute $y$ inexactly, i.e., $y$ is obtained by a floating-point LP solve, using the standard error tolerances. This suffices, since slightly negative dual multipliers can be set to zero and any conic combination of the constraints is globally valid.

To aggregate the constraints we also use floating-point arithmetic. We start from the constraint $\hat{a}^\T x\leq \hat{b}$ for $\hat{a}=0$ and $\hat{b}=0$, which holds trivially. Now we add the constraints $y_j A_{\cdot j}^\T x\leq y_j b_j$ (where $A_{\cdot j}$ is the $j$th column of $A$) to this constraint one by one: For each $j$, we set $\hat{b}\gets \roundup{\hat{b}+y_j b_j }$ and the components of $\hat{a}$ are updated in the following way:
\begin{itemize}
	\item If $u_i\leq 0$, then we set $\hat{a}_{i}\gets
		      \roundup{\hat{a}_i+y_{j}A_{ji}}$.
	\item Otherwise, if $\ell_i\geq  0$, then we round $\hat{a}_{i}\gets \rounddown{\hat{a}_i+y_{j}A_{ji}}$.
	\item Otherwise, if $u_i\leq  \infty$, then we set
	      $\hat{b}\gets \roundup{\hat{b}+ \roundup{(\roundup{\hat{a}_i + y_{j}A_{ji} } - \rounddown{(\hat{a}_i + y_{j}A_{ji})})}u_i}$ and set $\hat{a}_{i} \gets  \roundup{\hat{a}_i + y_{j}A_{ji} }$.
	\item Otherwise, if $\ell_i\geq  -\infty$, then we set
	      $\hat{b}\gets \roundup{\hat{b}- \roundup{(\roundup{\hat{a}_i + y_{j}A_{ji} } - \rounddown{(\hat{a}_i + y_{j}A_{ji})})}\ell_i}$ and set $\hat{a}_{i} \gets  \rounddown{\hat{a}_i + y_{j}A_{ji} }$.
\end{itemize}

Note that this way of rounding guarantees that the constraint stays valid, and yields an approximation of the exact conflict constraint.
This approach has been used before to implement Gomory mixed integer cuts~\cite{EiflerGleixner2023}.

\subsection{Producing Certificates of Correctness}
MIP solvers are complex pieces of software, that might contain bugs.
Hence, it can be desirable to verify the correctness of a solution that was found using a MIP solver, especially when finding exact solutions is important. While checking feasibility is easy, verifying optimality is harder, since this requires checking many nodes in the branch-and-bound tree without using the solver. To make it possible to verify optimality, a MIP solver can be extended to emit a certificate of optimality \cite{cheung_verifying_2017}. The certificate contains all derived bounds, plus a reason for why each bound holds. These certificates can be verified using a simple piece of software that does not make use of the solver.

SCIP has been extended to generate certificates that can be verified using the program VIPR~\cite{cheung_verifying_2017}. A VIPR certificate contains all the initial problem constraints and any derived constraints.  Each derived constraint contains a reason, that justifies why the constraint holds. Previously derived  constraints can be referred to by their line number in the certificate. In VIPR certificates there are two deduction rules to derive an inequality:

\begin{itemize}
	\item {\bfseries Linear combination.} If a constraint is a linear combination of previously derived constraints, it has to be valid. To certify the derivation, the line numbers of these constraints are printed, along with the corresponding coefficients.
	\item {\bfseries Rounding.} If $c^\T x \leq b$ has been previously derived such that $c$ is zero for all indices that correspond to non-integral variables, then $\sum_{i}\lfloor c_i \rfloor x_i$ must be integral for any feasible solution $x$. Hence, $\sum_{i}\lfloor c_i \rfloor x_i \leq \lfloor b \rfloor$ holds.  To certify the derivation, the line number of the constraint  $c^\T x \leq b$ is printed to the certificate. This procedure is called a \emph{Chvátal-Gomory cut} in the context of pure integer programming.
\end{itemize}
Additionally, there are deduction rules that allow to encode a branching proof.

To be able to verify bounds that were derived using constraint propagation or dual proof analysis, these bounds need to be written to the certificate in terms of the above derivation rules.
Because \cref{eq:act_bound1,eq:act_bound2} are both linear combinations of valid variable bounds, constraint propagation steps can be added to the certificate using the `linear combination' deduction rule. If $a_k>0$ the constraint is
\begin{align}
	x_k \geq \frac{1}{a_k}\left(b_j-\sum_{i\in [n] \setminus \{k\}}^n\max(a_i\ell_i, a_iu_i)\right).
\end{align}
For each $i\neq k$, the line number of $x_i\geq \ell_i$ is printed to the certificate if $a_i\ell_i> a_iu_i$. Otherwise, the line number of $x_i\leq u_i$ is printed. Finally, the line number corresponding to the constraint $\sum_{i=1}^n a_ix_i \geq b_j$ is printed as well. For each of the constraints, the corresponding coefficient is $1/a_k$.

If $x_k$ is an integer variable and for all $i$ with nonzero $a_i$, both $a_i$ and $x_i$ are integral, a rounding derivation is printed to the certificate, certifying
\begin{align}
	x_k \geq \big\lceil \frac{1}{a_k}(b_j-\sum_{i\in [n] \setminus \{k\}}^n\max(a_i\ell_i, a_iu_i))\big\rceil.
\end{align}

In dual proof analysis, only linear combinations of globally valid constraints are added as constraints. However, these constraints are computed using safe rounding methods instead of exact arithmetic. For that reason, the derived constraints can not be computed directly using the given multipliers and the ``linear compbination'' deduction rule.

Instead, we let the solver write just the linear combination corresponding to the original constraint to the certificate along with a list of the current bounds for all variables appearing in the constraint. Then afterwards a post-processing tool, called \texttt{viprcomplete} is applied to the certificate to repair the linear combination of all of the constraints that have slight inaccuracies due to the safe rounding procedure. One of the advantages of this approach is that the repair step only needs to be done for constraints that will actually be used in the solving process. For details, we refer to \cite{EiflerGleixner2023}.

\newcommand{\fpdef}{fp-baseline}
\newcommand{\fpprop}{fp-prop}
\newcommand{\fpconf}{fp-prop+conflict}

\newcommand{\exdef}{baseline}

\newcommand{\exprop}{Propagation}
\newcommand{\expropfreq}{Propagation (every 3 nodes)}
\newcommand{\expropcont}{Propagation (every 3 nodes) + continuous variables}

\newcommand{\exconf}{prop+conflict}
\newcommand{\exconffreq}{prop-freq3+conflict}
\newcommand{\exconfnosb}{prop+conf-nosb}

\section{Computational Study}
\FloatBarrier
\label{sec:comp}
To analyze the actual impact of constraint propagation and dual proof analysis in the context of exact MIP solving,
we implemented these techniques in the exact variant of SCIP and compared them with the impact of the same techniques in the floating-point version of SCIP. Our runtime experiments do not include the time for certificate generation in order to be as comparable as possible with the floating-point version. We stress that while certificate generation does incur a computational cost, it does not change the solving path and therefore does not affect the solvability of instances. We refer to~\cite{eifler_computational_2022,EiflerGleixner2023} for a detailed discussion of the computational cost of certificate generation.

\subsection{Experimental Setup and Test Set}
The experiments were all performed on a cluster of Intel Xeon Gold 5122 CPUs with
3.6 GHz and 96 GB main memory. For all symbolic computations, we use
the GNU Multiple Precision Library (GMP) 6.1.2 \cite{GMP}.
All compared algorithms are implemented within \scipv~\cite{scip8zen}, using \soplexv~\cite{soplex6zen} as both the floating-point and the exact LP solver. For presolving in exact arithmetic, we use \papilov\cite{papilo20zen} and disable all other presolving steps.
Our proposed algorithms are freely available on GitHub.\footnote{As part of the development version of exact \scip mirrored under \url{https://github.com/scipopt/scip/tree/exact-rational}.}

We use the \miplib~2017 benchmark test set~\cite{MIPLIB}; in order to save computational effort, we exclude all those that could not be solved by the floating-point default version of \scipv within two hours. For the remaining $132$ instances we use three random seeds, making the size of our test set $396$.
We report aggregated times in shifted geometric mean with a shift of $1$ second, and node numbers in shifted geometric mean with a shift of $100$ nodes. For all tests, a time limit of two hours is used.

\subsection{Experimental Results}
\begin{figure}[t]
	\centering
\includegraphics[width=.9\textwidth]{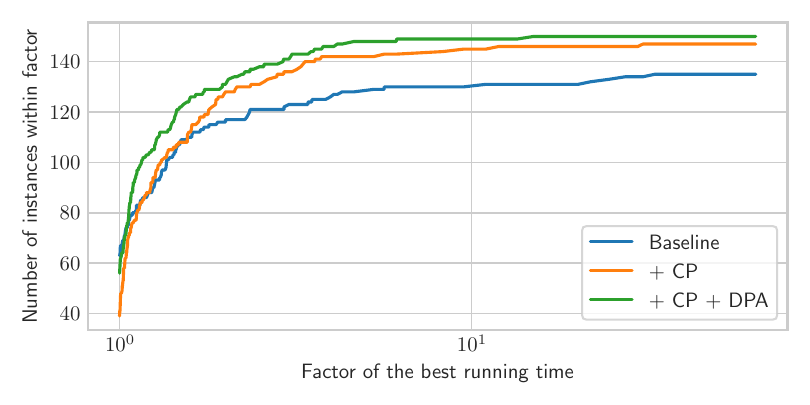}

		\caption{Performance profile for different configurations. For each configuration it is shown how many instances achieved a running time within a given factor of the shortest running time for that instances. Only the instances that could be solved by at least one of the given configurations are included in these statistics. \label{fig:performanceprofile}}
      \bigskip
\end{figure}

The experimental results can be found in \Cref{table:results}.
From the table we see that enabling propagation leads to an improvement in the running time of 11.6\%.
Also, 12 more instances are solved to optimality within the time limit.
The number of nodes decreases by 13.6\%.
Enabling dual proof analysis decreases the running time by another 13.4\% and the node count by 17.7\%.
3 more instances are solved to optimality within the time limit.
Together this comes down to a 23.4\% decrease in running time, a 28.9\% decrease in node count and 15 more instances that can be solved. From the performance profile in \cref{fig:performanceprofile} it is also clear that enabling constraint propagation and conflict analysis consistently speeds up the solver.
We also measured the total time that is spent in propagation. As can be seen from the table, this amount is small in comparison with the total time.

To compare this with the impact that these techniques have in the floating-point setting, we ran the same experiments on floating-point SCIP.
In these experiments we disabled the features not present in the exact version of SCIP: all cutting planes, graph-based conflict analysis, restarts, presolving (except for PAPILO), custom propagation rules.
As shown in \Cref{table:results_fp} in the floating-point setting, enabling propagation reduces solving time by 12.7\%.
The node count increases by 6.2\% and 11 more instances are solved to optimality. Dual proof analysis reduces the running time by another 36.5\% and the node count by 67.1\%. 20 more instances are solved to optimality.
Together this comes down to a decrease in the running time of 44.5\%, a 65.1\% decrease in node count and 31 more instances being solved to optimality.

\begin{table}[t]
   \centering
\begin{tabular*}{0.85\textwidth}{@{\extracolsep{\fill}}lrrlrrrl}
      \toprule
      {} & \# Solved & \multicolumn{4}{c}{Time} & Nodes & (rel) \\
      \cmidrule{3-6}
      Settings   &           & Total &(rel) & CP & DPA & \\
      \midrule
      Baseline   & 135       & 759.21 & --- & ---     & ---      & 8735.1 & --- \\ + CP       & 147       & 671.22 &(0.88)& 9.11    & ---      & 7550.4 &(0.86)\\ + CP + DPA & 150       & 581.46 &(0.77)& 9.22    & 5.44     & 6211.1 &(0.71)\\ \bottomrule
   \end{tabular*}
   \caption{Experimental results comparing the performance of exact SCIP with and without constraint propagation (CP) and dual proof analysis (DPA) enabled. All times shown are in seconds. For the times and node counts, the shifted geometric with shift 1 seconds or 100 nodes respectively is used. Columns (rel) show times and nodes relative to the baseline. Only the instances that could be solved by at least one of the given configurations are included in these statistics.	\label{table:results}}
   \bigskip
\end{table}
 
\begin{table}[t]
	\centering
\begin{tabular*}{0.8\textwidth}{l@{\extracolsep{\fill}}rrlrl}
		\toprule		Settings & \# Solved &  Time & (rel) &  Nodes & (rel) \\
\midrule
		Baseline      & 164   &   584.84 &  --- &       11415.0 & ---\\
		+ CP          & 175   &   510.85 &  (0.87)  &     12119.4 & (1.06) \\
		+ CP + DPA    & 195   &   324.41 &  (0.55)  &      6986.1 & (0.61) \\
		\bottomrule
	\end{tabular*}
	\caption{Experimental results comparing the performance of floating-point SCIP with and without constraint propagation (CP) and dual proof analysis (DPA) enabled. All times shown are in seconds. For the times and node counts, the shifted geometric with shift 1 seconds or 100 nodes respectively is used. Columns (rel) show times and nodes relative to the baseline.
		Only the instances that could be solved by at least one of the given configurations are included in these statistics.
		\label{table:results_fp}}
      \bigskip
\end{table}
 
It is clear that the performance boost is significantly larger in the floating-point setting. The same is true for the increase in the number of instances solved. While the speedup from only enabling propagation is comparable, we observe a much greater reduction in the number of nodes, and a much larger speedup from conflict analysis in the floating-point setting. We see several reasons contributing to this difference.

First, performance variability plays a part and the fact that we look at different subsets of the instances.
Second, propagation takes up more time in the exact solving mode. Although the shifted geometric mean of the propagation times shown in \Cref{table:results} is not large, there exist instances where exact propagation even takes up a major portion of the solving time.
This is the case for instances in which propagation leads to a large number of bound changes, since applying these bound changes is computationally more expensive in the exact setting. This is because, while the propagation procedure itself is implemented using floating-point arithmetic, applying a bound change in exact SCIP is currently done using exact arithmetic. This means every floating-point bound change needs to be converted to a rational number first, and then that rational number needs to be used to update the exact bound.
This can slow down the process significantly. If we only look at the solving time without the time spent in propagation, the speedup from enabling propagation is $17.1\%$, and the speedup from enabling conflict analysis is $28.3\%$.

There are more reasons that might lead to a lower impact of dual proof analysis in the exact setting.
The most important one is that constraint propagation is performed using floating-point arithmetic, and that we cannot determine infeasibility within tolerances in the exact setting.
This might sometimes lead to derived bounds not being precise enough to be useful in the exact settings or conflicts not being able to prove infeasibility.
For example, it is possible that in floating-point mode such an inexact bound can be used to decide the infeasibility of a node, whereas in exact mode the corresponding LP still needs to be solved. In numerically difficult cases, this might even lead to incorrect cutoffs in the floating-point setting.
It is clear that this is a trade-off that we can never fully avoid, guaranteeing correctness while sacrificing performance.

Another possible reason that the techniques are slightly less effective in exact SCIP could be that the derived variable bounds on non-integral variables can have a large denominator.  This can slow down solving LPs exactly, as has been observed for safe cutting plane generation in~\cite{EiflerGleixner2023}.
We experimented with two variants that try to prevent this from happening:
\begin{enumerate}
	\item Limit the size of the denominators of the derived variable bounds for continuous variables to a fixed upper bound.
	\item Disable constraint propagation for continuous variables.
\end{enumerate}
However, we did not observe any positive effects on the solving times. This does not mean that the described negative effects do not occur, but rather that all in all, by enabling propagation for non-integral variables and allowing arbitrary denominators in their bounds, we gain more than we lose.

\section{Conclusion and Outlook}
\label{sec:conclusion}

In this study, we investigated the feasibility and impact of constraint propagation and dual proof analysis in the exact MIP solver SCIP.
We found that enabling both techniques decreases the running time by 23.4\% and the number of nodes by 28.9\%.
This performance boost is significant, but less than the 44.5\% decrease in running time and 65.1\% decrease in node count that we observed in the equivalent floating-point setting. Such weaker performance is partially the price to be paid for exactness, coming from weaker inequalities due to directed rounding and the fact that we cannot determine infeasibility within tolerances in the exact setting.

Still, we see several future research opportunities to further improve these results. On the implementation side,
we believe there is potential to decrease the time spent in propagation by adding support for floating-point bound changes in exact SCIP.
Also, due to technical reasons, it is currently not possible to apply propagation within strong branching in the exact solving mode. Enabling this could also be beneficial.

Algorithmically, incorporating the strengthening techniques from \cite{witzig_computational_2021} for conflict constraints also in the exact setting is likely to yield additional performance improvements. 
Finally, implementing graph-based conflict analysis \cite{achterberg_conflict_2007} would be straightforward due to its combinatorial nature and is sure to yield positive impact in the exact setting.
More open-ended is the question if sporadically recomputing the activities also in the exact setting could lead to tighter bounds and thus better performance.

\newcommand{\etalchar}[1]{$^{#1}$}

\end{document}